\documentclass[12pt, reqno]{amsart}

\usepackage{graphics}

\newtheorem{theorem}{Theorem}
\newtheorem{proposition}[theorem]{Proposition}
\newtheorem{lemma}[theorem]{Lemma}
\newtheorem{corollary}[theorem]{Corollary}

\theoremstyle{definition}
\newtheorem{definition}[theorem]{Definition}

\theoremstyle{remark}
\newtheorem{remark}[theorem]{Remark}

\numberwithin{equation}{section}

\newcommand{\Z}{{\mathbf Z}}

\newcommand{\R}{{\mathbf R}}
\newcommand{\C}{{\mathbf C}}

\newcommand{\kk}{{\mathbf k}}
\newcommand{\cat}{{\rm {cat }}}

\newcommand{\I}{{\mathcal I}}
\newcommand{\K}{{\mathcal K}}

\begin{document}

\title{Topology of billiard problems, I}         
\author{Michael Farber}        
\address{Department of Mathematics, Tel Aviv University, Tel Aviv, 69978, Israel}
\email{farber@math.tau.ac.il}

\date{May 15, 2000}          

\thanks{Partially supported by the US - Israel Binational Science Foundation and by the Minkowski Center for Geometry. 
A part of this work was done while the author was visiting IHES}

\maketitle
\begin{abstract} Let $T\subset \R^{m+1}$ be a strictly convex domain bounded by a smooth 
hypersurface $X=\partial T$. 
In this paper we find lower bounds on the number of billiard trajectories in $T$, which have a prescribed 
initial point $A\in X$, a prescribed final point $B\in X$ and make
a prescribed number $n$ of reflections at the boundary $X$. 
We apply a topological approach based on calculation of cohomology rings
of certain configurations spaces of $S^m$.

\end{abstract}

\section{\bf Introduction}      

In the early 1900s, G.D. Birkhoff initiated the mathematical theory of convex plane billiards.
His main interest was in estimating the number of periodic billiard trajectories. 
He has pioneered the use of topological methods, based on the variational reduction and using
the critical point theory. 

Periodic trajectories in convex billiards in Euclidean spaces of dimension $>2$  were studied in \cite{Ba} and \cite{FT}.
The high-dimensional problem also allows approach based on the critical point theory, and the main difficulty lies in 
more complicated topology of the appropriate configuration space. Thus, the major effort 
of \cite{FT} was in computing 
the cohomology algebra of the cyclic configuration space of the sphere.

The purpose of this paper is to strengthen the estimates obtained in \cite{FT}.
In particular, we will obtain linear in $n$ estimates on the number of $n$-periodic billiard trajectories,   
improving the logarithmic estimates of \cite{FT}. 
On the other hand, in this work we will study a larger variety of billiard
problems: besides the periodic trajectories we will be also interested in the number of ways
the billiard ball can be brought from a given initial point to a given final point after making a prescribed 
number $n$ of reflections at the boundary of the billiard domain.

Let $X\subset \R^{m+1}$ be a closed smooth strictly convex hypersurface. 
The $(m+1)$-dimensional convex body $T$, bounded by $X$ will serve as our billiard table. 
The billiard ball is a point which moves in $T$ in a straight line, except when it hits $X= \partial T$,
where it rebounds making the angle of incidence equal the angle of reflection.
In other words, if $A, B, C\in X$ are three subsequent reflection points then the normal to $X$ at point $B$ bisects
the angle between the vectors $BA$ and $BC$.

We will consider the following billiard problems. 

{\bf Problem A:} {\it Given two distinct points $A,B \in X$ and a number $n$, estimate the number of billiard trajectories
inside $X$ which start at point A, end at point $B$ and make $n$ reflections at the hypersurface $X$.}

{\bf Problem B:} {\it Estimate the number of billiard trajectories inside $X$, which start and end at a given point $A\in X$ and
make a prescribed number $n$ of reflections at the hypersurface $X$.} 

Problem B deals with closed billiard trajectories.
It is clear that any closed billiard trajectory starting and ending at $A\in X$ determines another closed billiard trajectory
which is obtained by passing the same route in the reverse order. This explains that there is a natural 
$\Z_2$-action on the set
of closed billiard trajectories, and in Problem B one actually asks about {\it the number of 
$\Z_2$-orbits} of closed billiard trajectories. 

Using this $\Z_2$-symmetry we will give a better estimate for 
Problem B than the estimate for Problem A given by Theorem \ref{thm1} below . 

{\bf Problem C:}  {\it Estimate the number of $n$-periodic billiard trajectories inside the billiard domain $T$.}

In \cite{FT} we showed that the number of $n$-periodic billiard trajectories is at least 
$[\log_2 (n-1)]+m$, and it is at least $(n-1)m$ in the generic case. 

\centerline{\includegraphics{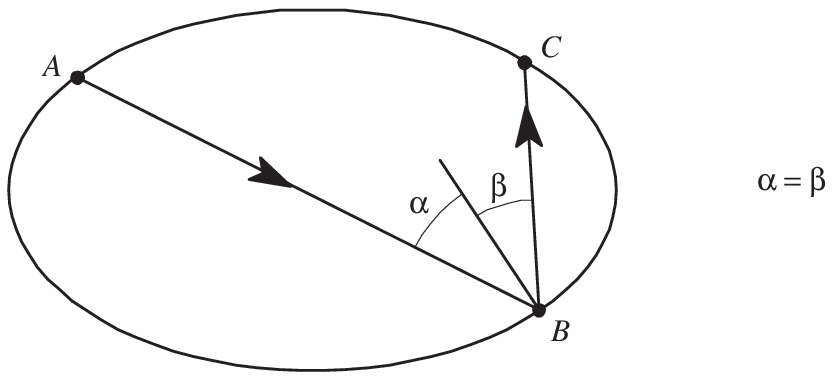}}

In Part I of this paper we will give an answer to  Problem A. 

Our purpose is to prove the following Theorem.

\begin{theorem}\label{thm1}
Let $X\subset \R^{m+1}$ be a closed smooth strictly convex hypersurface, $A, B\in X$ two distinct points.
Then for any integer $n$ the number of billiard trajectories inside $X$, which start at $A$, 
end at $B$ and make
$n$ reflections is at least $n+1$, if $m$ is odd, and $[(n+1)/2]+1$, if $m$ is even.
 In a generic situation (cf. Definition \ref{generic} below) for any $m$ and $n$  
the number of billiard trajectories inside 
$X$, which start at $A$, end at $B$ and make
$n$ reflection, is at least $n+1$.
\end{theorem}

Here $[x ]$ denotes the integer part of $x$, i.e. the largest integer not exceeding $x$.

Consider the following simple example. Let  $X=S^{m}\subset \R^{m+1}$ be the unit sphere. 
Any billiard trajectory $A=A_0, A_1, \dots, A_n, A_{n+1}=B$ must lie in a 2-plane passing through the center of the sphere $O$.
If the endpoints $A, B$ are distinct and
not antipodal then there is a unique 2-plane passing through $A, B$ and $O$; the circle $L$,  
the intersection  of this 2-plane with $S^m$, must contain all the reflection points $A_1, \dots, A_n$.
Fix an orientation on $L$ and let $\phi\in (0, 2\pi)$ be the angle from $A$ to $B$. Then the angle between $A_i$ and $A_{i+1}$
must be independent of $i$ and may take the values
\[\alpha_k = (\phi +2\pi k)/({n+1}), \quad\mbox{where}\quad k=0, 1, \dots n.\]
Hence, we see that in this example 
there exist precisely $n+1$ billiard trajectories starting at $A$, ending at $B$ and making $n$ reflections. 

This example shows that the statement of Theorem \ref{thm1} for the generic case and for 
$m$ odd cannot be improved. It looks reasonable to conjecture that for even $m$ the lower bound on the number
of billiard trajectories is also $n+1$.

Problems B and C will be studied in Part II.

I would like to thank S. Tabachnikov for useful discussions.

\section{\bf Billiard ball problem and Lusternik - Schnirelman category of configuration spaces}

In this section we will use the variational method of G. Birkhoff
to show that the problem of estimating the number of billiard trajectories can be
reduced to a topological problem of estimating the Lusternik-Schnirelman category of a space of
configurations of $n$ points on the sphere $S^m$.
 
Let $X$ be a manifold. Suppose that $A, B\in X$ are two fixed points.
The symbol $G(X;A,B,n)$ will denote  the subspace of 
the Cartesian power $X^{\times n}= X\times X\times \dots \times X$, consisting of 
the configurations $(x_1, \dots, x_n)\in X^{\times n}$, such that $x_i\ne x_{i+1}$ for all $i= 1, \dots, n-1$, 
and  $A\ne x_1$ and $x_{n} \ne B$. 
In the case $A\ne B$ we will call $G(X;A,B,n)$ {\it the open string configuration space}.
The space $G(X;A,A,n)$ is {\it the closed string configuration space}. 

The configuration space $G(X;A,B,n)$ is closely related to the {\it cyclic
configuration space}  $G(X,n)$ introduced in \cite{FT}, which consists of all 
 $n$-tuples $(x_1,  \dots, x_n)$ of points of $X$ satisfying:
\begin{eqnarray*}
x_i\ne x_{i+1}\quad {\rm for}\quad i=1, 2, \dots, n-1,\quad \mbox{and}\quad x_n\ne x_1.\label{cyclic}
\end{eqnarray*}

Let $X \subset {\R}^{m+1}$ be a smooth closed strictly
convex hypersurface, the boundary
of the billiard table. Denote by
\begin{eqnarray*}
L_X: G(X;A,B,n)\to \R
\end{eqnarray*}\label{length}
the perimeter length function, taken with the minus sign,
\begin{eqnarray*}
L_X(x_1,\dots,x_n) = -\sum_{i=0}^n |x_i - x_{i+1}|, \quad i=0, \dots,n,
\end{eqnarray*}
where $(x_1, x_2, \dots, x_n)\in G(X;A,B,n)$ and the distance $|x_i-x_{i+1}|$
is measured in the ambient Euclidean space $\R^{m+1}$. Here we understand that $x_0=A$ and $x_{n+1}=B$.
The function $L_X$ is smooth. The reason for the minus sign will become clear later. 

The following Lemma is well known.

\begin{lemma} A point $(x_1, \dots, x_n)\in G(X;A,B,n)$ is a critical point of $L_X$
if and only if the sequence 
$A, x_1, \dots, x_n, B$ determines a billiard trajectory inside $X$,
starting at point $A$ and ending at point $B$.
\end{lemma}
\begin{proof}
An easy calculation shows that a configuration $(x_1, \dots, x_n)\in G(X;A,B,n)$ is a critical point of $L_X$
if and only if for any $i=1, 2, \dots, n$ the vector
\[\frac{x_i-x_{i-1}}{|x_i-x_{i+1}|} + \frac{x_i-x_{i+1}}{|x_i-x_{i+1}|}\]
is orthogonal to the tangent space $T_{x_i}(X)$. The last condition is clearly equivalent to the requirement that
the normal to $X$ at $x_i$ bisects the angle between $x_ix_{i-1}$ and $x_ix_{i+1}$.
\end{proof}

\begin{definition}\label{generic} The data $(X,A,B,n)$ will be called generic, if the associated perimeter length function 
$L_X: G(X;A,B,n)\to \R$ has only Morse critical points.
\end{definition}

Compare \S 4 of \cite{FT}. 

Since $X$ is homeomorphic to $S^m$, the space $G(X;A,B,n)$ is homeomorphic to $G(S^m;A,B, n)$.
The
shape of the billiard
domain $X$ becomes encoded in the function $L_X: G(S^m;A,B,n)\to \R$, and the
problem of estimating the number of billiard trajectories inside $X$, which start at $A$ and end at $B$,
turns into a problem of
Morse-Lusternik-Schnirelman theory.
The difficulty is that we cannot apply
the Morse-Lusternik-Schnirelman
theory directly to $G(X;A,B,n)$ since
this manifold is not compact.

To avoid this difficulty, we will replace
$G(X;A,B,n)$ by a compact manifold with boundary $G_{\varepsilon} (X;A,B,n) \subset
G(X;A,B,n)$, where $\varepsilon >0$ is small enough and
\begin{eqnarray}
G_{\varepsilon} (X;A,B,n) = \{(x_1, \dots, x_n) \in X^{\times n} :
\prod_{i=0}^n|x_i - x_{i+1}| \geq \varepsilon \};
\end{eqnarray}
here $x_0=A$ and $x_{n+1}=B$.
A similar approach can be found in \cite {Ba} and in \cite {CS}, \cite {KT}
for the two-dimensional case; cf. also \cite{FT} for the periodic case.

\begin{proposition} \label{prop2} 
If $\varepsilon>0$ is sufficiently small
then:
\begin{enumerate}
\item[(a)] $G_{\varepsilon} (X;A,B,n)$ is a smooth manifold with boundary;
\item[(b)] the inclusion $G_{\varepsilon} (X;A,B,n)\subset G(X;A,B, n)$ is a homotopy equivalence;
\item[(c)] all
critical points of $L_X: G(X;A,B, n)\to \R$ are contained in $G_{\varepsilon}
(X;A,B,n)$;
\item[(d)] at every point of  $\partial G_{\varepsilon} (X;A,B,n)$,
the gradient of $L_X$ has the outward direction.
\end{enumerate}
\end{proposition}
This statement is analogous to Proposition 4.1 from \cite{FT}. The proof given in
\cite{FT} also applies in this case. The only modification is that in the case $A\ne B$, 
 the arguments of the proof
of Proposition 4.1 of \cite{FT}, which explain that a critical configuration cannot lie entirely in a small 
neighborhood of $X$, become redundant.

\begin{corollary}\label{cor1}
Let $X\subset \R^{m+1}$ be a smooth strictly convex hypersurface, $A, B\in X$ two fixed points.
For any $n\ge 0$, the number of billiard trajectories inside $X$, which start at $A$, end at $B$ and make $n$ reflections is at least
$\cat(G(S^m;A,B,n)),$
the Lusternik - Schnirelman category of the open string configuration space of the sphere $S^m$.
\end{corollary}

\begin{proof} Choose $\varepsilon >0$ small enough such that the conclusions of Proposition \ref{prop2} hold.
Since at the points of the boundary $\partial G_{\varepsilon} (X;A,B,n)$
the gradient of $L_X$ has the outward direction, the critical point theory  for manifolds with boundary \cite{BM}
applies; the conclusion is that the critical points of the restriction 
$L_X|_{\partial G_{\varepsilon} (X;A,B,n)}$ should be 
ignored, and the 
number of critical pints of $L_X$ lying in the interior of $G_{\varepsilon}(X;A,B,n)$ is at least 
the category of $G_{\varepsilon}(X;A,B,n)$.  Since 
$\cat G_{\varepsilon}(X;A,B,n) =  \cat G(S^m;A,B,n)$ (because of statement (b) of Proposition \ref{prop2}),
the number of billiard trajectories inside $X$, which start at $A$, end at $B$ and make $n$ reflections is at least
$\cat(G(S^m;A,B,n)).$
\end{proof}

In the closed case, i.e. assuming that the endpoints are equal
$A=B$, we may use $\Z_2$-symmetry to give a better estimate.
This result will appear in Part II.

\section{\bf Spectral sequence computing cohomology of the open string configuration space of a manifold }

The following Theorem yields a spectral sequence computing the cohomology algebra of the open string configuration space
$G(X;A,B,n)$, where $X$ is an arbitrary manifold.
It is a Leray spectral sequence associated to the embedding $G(X;A,B,n)\to X^{\times n}= X\times X\times \dots \times X$
(the $n$-th Cartesian power). 

This method was first suggested by B. Totaro \cite{To} for the usual configuration space
(i.e. for the space of all configurations $(x_1, \dots, x_n)\in X^{\times n}$ with $x_i\ne x_j$ for all $i, j$). In \cite{FT}
we used a similar spectral sequence for the cyclic configuration space $G(S^m,n)$.
 
$\kk$ will denote a field.

\begin{theorem} \label{openstring1}
Let $X$ be a connected oriented manifold of dimension $m>1$ and let $A, B\in X$ be two distinct points.

(A) There exists a spectral sequence of bigraded differential algebras, which converges to $H^\ast(G(X;A,B,n);\kk)$,
whose $E_2$-term is the quotient of the bigraded commutative algebra
\[H^\ast(X^{\times n};\kk)[s_0, s_1, \dots, s_n],\]
where $H^p(X^{\times n};\kk)$ has bidegree $(p,0)$ and each generator $s_i$ has bidegree $(0,m-1)$, by the relations
$$
\begin{array}{ll}
s_i^2 =0, & \mbox{for}\quad i=0, 1, \dots, n,\\
s_0s_1\dots s_n =0,&\\
p_1^\ast(v) s_0=0, &\\
p^\ast_i(v)s_i =p^\ast_{i+1}(v)s_i, \quad &\mbox{for}\quad i=1, 2, \dots, n-1, \\
p_n^\ast(v) s_n =0,&
\end{array}
$$
where $v\in H^\ast(X;\kk)$ denotes an arbitrary cohomology class of positive degree and $p_j: X^{\times n}\to X$ denotes the
projection onto the $j$-th factor, $j=1, 2, \dots, n$. 

(B) The first nontrivial differential is $d_m$, where $m=\dim X$. It acts
by 
$$
\begin{array}{l}
d_m(s_0) = (-1)^mp_1^\ast([X]), \\
d_m(s_i) = q_i^\ast(\Delta), \quad i=1, 2, \dots, n-1,\\
d_m(s_n) = p_n^\ast([X]),\\
d_m(H^\ast(X^{\times n};\kk))=0,
\end{array}
$$
where $q_j: X^{\times n} \to X\times X$ denotes the projection onto the factors $j$ and $j+1$, 
$[X]\in H^m(X;\kk)$ is the fundamental
class, and
$\Delta\in H^m(X\times X;\kk)$ denotes the cohomology class of the diagonal.
\end{theorem}

\begin{proof}[Proof of Theorem \ref{openstring1}] 
Consider the inclusion $\psi: G(X;A,B,n) \to X^{\times n}$ and the Leray
spectral sequence \cite{Bor} of the continuous map $\psi$
$$E_2^{p,q}\, =\, H^p(X^{\times n}; R^q \psi_{\ast} {\kk}) \Rightarrow
H^{p+q} (G(X;A,B,n);\kk),$$
where $R^q \psi_{\ast} {\kk}$ is the sheaf on $X^{\times n}$ associated
with the presheaf
$$U \mapsto H^q(U \cap G(X;A,B,n);\kk).$$

To describe the sheaves $R^q \psi_{\ast} {\kk}$, consider partitions of the set $\{0, 1,\dots,n, n+1\}$ into
intervals, that is, subsets of the form $\{i,i+1,i+2, \dots, i+j\}$.
For any such partition $J$ we denote by $X_J$ the subset of $X^{\times n}$,
consisting
of all configurations $c=(x_1, x_2, \dots, x_n)\in X^{\times n}$,
satisfying the conditions: 
$$\begin{array}{l}
x_i = x_j, \quad\mbox{ if $i$ and $j$ lie in the same interval of the partition $J$};\\
x_i=A,\quad\mbox{ if index $i$ lies in the same interval with $0$};\\
x_i=B,\quad\mbox{ if index $i$ lies in the same interval with $n+1$}.
\end{array}$$
Given two interval partitions $I$ and $J$, we  say that $J$ {\it
refines} $I$ and write $I\prec J$
if the intervals of $I$ are unions of the intervals of $J$.
We denote by $|J|$ the number of intervals in the partition $J$.
Note that $I\prec J$ implies $X_I\subset X_J$ and $|I|\le |J|$.
For the partition $J$ with $|J|=1$ holds $X_J =\emptyset$ (since we assume that $A\ne B$).
If $|J|=2$ then $X_J$ is a single point. For $|J|>2$ the space
$X_J$ is homeomorphic to the Cartesian power $X^{\times (|J|-2)}$.

As in \cite{FT}, we will denote by $D(X,n)$  the subset of $X^{\times n}$ satisfying the conditions
$x_i \neq x_{i+1}$ for $i=1, \dots, n-1$. 
The configuration space 
$D({\R}^m,n)$ is homotopy equivalent to the product of spheres
$(S^{m-1})^{\times (n-1)}$. A homotopy equivalence $D({\R}^m,n)\to (S^{m-1})^{\times (n-1)}$ 
is given by the map
\begin{eqnarray}\label{canonical}
(x_1, \dots, x_n)\mapsto (\frac{x_1-x_2}{|x_1-x_2|}, \dots, \frac{x_{n-1}-x_{n}}{|x_{n-1}-x_{n}|}).
\end{eqnarray}
Fixing an orientation of the sphere $[S^{m-1}]\in H^{m-1}(S^{m-1};\kk)$, 
determines a {\it canonical top-dimensional class} in $H^{(n-1)(m-1)}(D(\R^m,n);\kk)$, which is the pull-back of the
product $[S^{m-1}]\times \dots\times [S^{m-1}]$ under (\ref{canonical}). 

If $A, B\in X$ two points, we denote by $G(X;A,\emptyset,n)$ the subspace 
of $D(X,n)$ consisting of configurations $(x_1, \dots, x_n)$ with
 $x_1\ne A$; similarly we denote by $G(X;\emptyset, B,n)\subset D(X,n)$ the subspace of configurations with $x_n \ne B$.

Let $J$ be a partition of $\{0, 1, 2, \dots, n+1\}$ on intervals of lengths
$j_1,...,j_r$, and let 
\begin{eqnarray*}
c =(x_1, x_2, \dots x_n)\, \in X_J, \quad c\notin \bigcup_{{I\prec J}, {I\ne J}}X_I.
\end{eqnarray*}
We claim that the stalk of the sheaf $R^q \psi_{\ast} {\kk}$ at $c$ equals
$$
(R^q \psi_{\ast} {\kk})_c =
H^q(D({\R}^m,j_1) \times \dots \times D({\R}^m,j_r);\kk).
$$
Indeed, by definition, this stalk is $H^q(U \cap G(X;A,B,n);\kk),$
 where $U$ is a small open ball around $c$. If $c=(x_1, x_2, \dots x_n)$
then we may choose points
$y_1=A, y_2, \dots, y_{r-1}, y_r=B\in X$, one for each interval of $J$, so that
$x_i=y_{j_s}$ if $i$ belongs to the
$s$-th interval.  Let $U_j\subset X$ be a small open neighborhood of $y_j$,
so that each $U_j$ is diffeomorphic to $\R^m$ and
the sets $U_j$ and $U_{j'}$ are disjoint when the points $y_j$
and $y_{j'}$ are distinct. Then we may take
$U=U_1^{\times j_1}\times U_2^{\times j_2}\times \dots \times U_r^{\times
j_r}$,
and our claim follows.

We see that $R^q \psi_{\ast} {\kk}$ vanishes unless $q$ is
a multiple of $m-1$
and 
$$
{\dim}\ (R^{s(m-1)} \psi_{\ast} {\kk})_c = 
\left\{
\begin{array}{l}
0, \quad \mbox{for}\quad s> n+2-|J|,\\ \\
 \left(
\begin{array}{c}
n+2-|J|\\ 
s
\end{array}
\right),\quad\mbox{for}\quad s\le n+2-|J|.
\end{array}
\right.
$$

For an interval partition $J$ of $\{0, 1, 2,\dots, n+1\}$ with $|J|>1$
 denote by $\varepsilon_J$ the constant sheaf
with stalk $\kk$ and support $X_J$. We claim: 

{\it for any $r=2, 3, \dots, n+2$, the sheaf
$R^{(n+2-r)(m-1)} \psi_{\ast} {\kk}$ is isomorphic to the direct sum of sheaves
\begin{eqnarray}
R^{(n+2-r)(m-1)} \psi_{\ast} {\kk}\, \simeq\, \bigoplus_{|J|=r}
\varepsilon_J,\label{split}\end{eqnarray}the sum taken over all interval
partitions $J$ with $
|J|=r$.}

To prove the claim, let
$I$ be an interval partition
of $\{0, 1,  \dots, n+1\}$ into intervals of length $i_1, i_2, \dots, i_s$,
where $s=|I|>1$. Then for any
interval partition $J$ into intervals of length $j_1, j_2, \dots, j_{r}$,
such that $I\prec J$, we have the
canonical inclusion
\[\nu_{JI}: D(\R^m,i_1)\times \dots \times D(\R^m,i_s) \to
D(\R^m, j_1)\times \dots \times D(\R^m, j_{r}).\]
The target space of map $\nu_{JI}$ has a canonical nonzero
$(n+2-r)(m-1)$-dimensional cohomology class (cf. above).
The induced map $\nu_{JI}^\ast$
on $(n+2-r)(m-1)$-dimensional cohomology with $\kk$ coefficients
is a monomorphism. Let $z_{JI}$ denote the image of the top-dimensional
canonical class under the induced map $\nu_{JI}^\ast$.
Then (similarly to Lemma 3 in \cite{To}) for a fixed $I$, the classes $\{z_{JI}\}$ form a linear basis of the
cohomology
$H^{(n+2-r)(m-1)}(D(\R^m,i_1)\times \dots\times D(\R^m,i_{s});\kk)$, where $J$ runs
over all partitions with
$I\prec J$  and $|J|=r$. 

Indeed, using the map into product of spheres (\ref{canonical}) we see that
a linear basis of the cohomology
$H^{(n+2-r)(m-1)}(D(\R^m,i_1)\times \dots\times D(\R^m,i_{s});\kk)$
form monomials $s_{a_1} s_{a_2}\dots s_{a_{n+2-r}}$
with $0\le a_1<a_2<\dots <a_{n+2-r}\le n+1$, such that for any $s=1, \dots, n+2-r$
the indices $a_s$ and $a_{s}+1$ belong to the same interval of partition $I$. Let $ J$ be the 
partition determined by the equivalence relation on $\{0,1, \dots, n+1\}$, where $a_s\sim a_s +1$.
Then $I\prec J$ and the above monomial coincides with the class $z_{JI}$.

Given a partition $J$ of $\{0, 1, \dots, n+1\}$ on intervals of 
length $j_1+1, j_2, \dots,j_{r-1}, j_{r}+1$, where $r>1$,
consider the commutative diagram
$$
\begin{array}{clc}
G(X;A,B,n) & \stackrel {\psi}{\longrightarrow} & X^{\times n}\\ \\
\downarrow & & \downarrow {\rm id}\\ \\
G(X;A, \emptyset , j_1)\times D(X,j_2)\times\dots\times G(X; \emptyset, B, j_{r}) & 
\stackrel {g_J}{\longrightarrow} & X^{\times n}
\end{array}
$$
formed by the natural inclusions. Define sheaf $\varepsilon'_J = R^{(n+2-r)(m-1)}{g_J}_\ast(\kk)$ over
$X^{\times n}$. 
We want to show that $\varepsilon'_J$ is isomorphic to $\varepsilon_J$,
i.e. it is the constant sheaf
with stalk $\kk$ and support $X_J$. First,
$\varepsilon'_J$
vanishes outside $X_J$ (since we are considering the cohomology of the top
dimension). Let
$U$ be a small open neighborhood of a point $c\in X_J\subset X^{\times
n}$, such that
$U=\prod U_i$, where all $U_i$ are small
open disks and $U_i=U_j$ if $i$ and $j$ lie in the same interval
of $J$. Then
$$\varepsilon'_J(U) = H^{s(m-1)}(D(U_{i_1}, j_1+1)\times D(U_{i_2},j_2)\times \dots \times
D(U_{i_{r}},j_{r}+1);\kk) \simeq \kk$$
(where $s=n+2-r$)
has a canonical element (cf. above). This gives a continuous section of $\varepsilon'_J$ over $X_J$, and
hence $\varepsilon'_J\simeq\varepsilon_J$.

The commutative diagram above gives a map of sheaves $\varepsilon_J\to
R^{(n+2-r)(m-1)}\psi_\ast(\kk)$,
and summing, we obtain a map of sheaves
\[\bigoplus_{|J|=r} \varepsilon_J\, \to\,  R^{q}\psi_\ast(\kk),\quad \mbox{where}\quad q=(n+2-r)(m-1),\]
which, as we have seen above, is an isomorphism on stalks;
hence it is an isomorphism, and the claim (\ref{split}) follows.

We arrive at the following description of the term $E_2$ of the Leray spectral
sequence
$$
E_2^{p, r(m-1)} =
\bigoplus_{|J|=n+2-r} H^p(X_J;\kk),
$$
where $J$ runs over all partitions of $\{0, 1, \dots, n+1\}$ with $|J|>1$.
In order to identify this description with the one given in the statement
of the theorem,
assign to a monomial $s_{i_1} \dots s_{i_r}$ with $0\le i_1< i_2<\dots <i_r\le n$
the equivalence relation on the set of indices $\{0, 1,\dots,n+1\}$ generated by
$$i_1 \sim i_1 +1, \, \, i_2 \sim i_2 + 1, \dots, i_r \sim i_r + 1.$$
This equivalence relation defines a partition $J$ of the set
$\{0, 1,\dots,n+1\}$ on $n+2-r$ intervals.
In view of the relations 
\begin{eqnarray}
\begin{array}{l}
p_1^\ast(v)s_0 = 0,\\
p_i^\ast(v)s_i = p_{i+1}^\ast(v)s_i, \quad\mbox{where}\quad i=1, \dots, n-1,\\
p_n^\ast(v)s_n = 0,
\end{array}
\end{eqnarray}
 the term
$H^p (X^{\times n};\kk) s_{i_1} \dots s_{i_r}$ is isomorphic to 
$H^p(X_J;\kk)$.

The monomial $s_0s_1\dots s_n$ corresponds to the partition $|J|=1$, which we should ignore since $X_J=\emptyset$; 
this explains the
relation $s_0s_1\dots s_n=0$.

Now we prove statement (B) concerning the differentials of the spectral sequence. 
The first nontrivial differential is $d_m$. To find $d_m$ it is enough to find the cohomology classes
$d_m(s_i)\in H^m(X^{\times n};\kk)$, where $i=0, 1, \dots, n$.

We will use functoriality of the Leray spectral sequence and the following well-known property. 
Let $Y$ be a manifold and let $Z\subset Y$ be a submanifold of codimension $m>1$ with oriented normal bundle. 
Consider the Leray spectral sequence
\[E_2^{p,q}\, =\, H^p(Y;R^q\phi_\ast\kk) \Rightarrow H^{p+q}(Y-Z;\kk),\]
of the inclusion $\phi: (Y-Z)\to Y$. The sheaf $R^{m-1}\phi_\ast\kk$ is the constant sheaf 
with support $Z$ and stalk $\kk$ for $q=m-1$ and it vanishes for all other values $q>0$. The only nonzero
differential
$d_m: E_2^{0,m-1} \to E_2^{m,0}$ acts as follows: the class $1\in H^0(Z;\kk) = E_2^{0,m-1}$ is mapped into
$d_m(1)=[Z]\, \in H^m(Y;\kk)$, the class dual to $Z$, where the same orientation of the normal bundle to $Z$ is used 
in order to trivialize the sheaf $R^{m-1}\phi_\ast\kk$ and to define the dual class $[Z]$. 

In order to show the first relation $d_m(s_0)=(-1)^mp_1^\ast([X])$, consider the diagram
$$
\begin{array}{ccc}
G(X;A,B,n) & \stackrel \subset \to & X^{\times n}\\
\downarrow & & \downarrow {\rm id}\\
(X-A)\times X^{\times (n-1)} & \stackrel \subset \to & X^{\times n}
\end{array}
$$
and apply the previous remark to the bottom row with $Y=X^{\times n}$ and $Z= A\times X^{\times (n-1)}$. 
The sign $(-1)^m$ appears as the degree of the antipodal map $S^{m-1} \to S^{m-1}$: the framing of the normal
bundle to $A\subset X$, which we use to define the fundamental class $[X]\in H^m(X;\kk)$, is antipodal to the framing
determined by (\ref{canonical}), which we use to trivialize the derived sheaf.

To obtain relations $d_m(s_i) = q_i^\ast(\Delta)$ with $i=1, \dots, n-1$ we use the commutative diagram
$$
\begin{array}{ccc}
G(X;A,B,n) & \stackrel \subset \longrightarrow & X^{\times n}\\
\downarrow & & \downarrow {\rm id}\\
X^{\times (i-1)}\times (X\times X - \Delta)\times X^{\times (n-i-1)} & \stackrel \subset \to & X^{\times n}
\end{array}
$$
and apply the remark above with $Y = X^{\times n}$ and 
$Z= X^{\times (i-1)}\times \Delta \times X^{\times (n-i-1)}$. The last relation $d_m(s_n)=p_n^\ast([X])$ follows similarly.

\end{proof}

We will apply Theorem \ref{openstring4} and Proposition 2.2 from \cite{FT}
to compute the integral cohomology of the configuration space $G(\R^m;A,B,n)$.
For any $i=0, 1, \dots, n$ we have the map
$$\phi_i: G(\R^m;A,B,n)\to S^{m-1}, \quad (x_1, \dots, x_n) \mapsto \frac{x_i-x_{i+1}}{|x_i-x_{i+1}|}\, 
\in S^{m-1},$$
where we understand $x_0=A$ and $x_{n+1}=B$.
Define the cohomology classes
$$s_i\in H^{m-1}(G(\R^m;A,B,n);\Z), \quad\mbox{as}\quad s_i =\phi_i^\ast([S^{m-1}]),\quad i=0, 1, \dots, n.$$

\begin{proposition}\label{cor13}
 For $m>1$ the algebra $H^\ast(G(\R^m;A,B,n);\Z)$ (where $A\ne B$)
is generated by cohomology classes
\[s_0, s_1, \dots, s_n\, \in H^{m-1}(G(\R^m;A,B,n);\Z)\]
and all relations between the classes $s_i$ are consequences of
$$
\begin{array}{l}
s_0^2 = s_1^2=\dots =s_n^2 =0,\quad s_is_j =(-1)^{m-1}s_js_i,\quad s_0s_1\dots s_n=0.\\
\end{array}
$$
\end{proposition}
\begin{proof} If we replace $\Z$ by a field $\kk$, the result follows directly from Theorem \ref{openstring4}.
In particular, we see that the dimension of the cohomology of $G(\R^m;A,B,n)$ 
do not depend on the field of coefficients. We conclude that
the integral cohomology of $G(\R^m;A,B,n)$ has no torsion and is nonzero only in dimensions divisible by $m-1$.

Consider the cyclic configuration space $G(\R^m, n+2)$ (cf. \cite{FT}) and the fibration
\begin{eqnarray}\label{fibration}
\quad \quad G(\R^m,n+2)\to G(\R^m,2)\simeq S^{m-1}, \quad (x_1, \dots, x_{n+2})\mapsto (x_{n+2}, x_{1}),
\end{eqnarray}
which has $G(\R^m;A,B,n)$ as the fiber. The nonzero rows of the Serre spectral sequence have numbers divisible by $m-1$;
also, the spectral sequence has only two columns $p=0$ and $p=m-1$. 
We obtain that all differentials of the spectral sequence vanish and the cohomology
of the fiber $H^\ast(G(\R^m;A,B,n);\Z)$ is isomorphic to the factor of the ring
$H^\ast(G(\R^m,n+2);\Z)$ with respect to the 
ideal generated by class $s_{n+2}$ (the pull-back of the fundamental class of the base). 
Comparing the above information with the structure of the ring $H^\ast(G(\R^m,n+2);\Z)$,
described in Proposition 2.2 of \cite{FT}, proves Proposition \ref{cor13}.
\end{proof}

\section{\bf Cohomology of open string configuration spaces of spheres}

In this section we state a theorem describing the cohomology 
of the configuration space $G(S^m;A,B,n)$, assuming that the points $A$ and $B$ are distinct.
We will see that the additive structure of the cohomology algebra $H^\ast(G(S^m;A,B,n);\kk)$ is similar for all $m$,
but the multiplication depends on the parity of the dimension $m$. Also, the
case $m=1$  is special, for the space $G(S^1;A,B,n)$ consists of $n+1$ path connected components and each is contractible
(cf. section \ref{proof}); in this case only zero-dimensional cohomology exists.

Let $\kk$ be a field.

\begin{theorem}\label{openstring4}
The cohomology $H^\ast(G(S^m;A,B,n);\kk)$ of the open string configuration space
 (where $A\ne B$) has additive generators 
\[\sigma_i \in \, H^{i(m-1)}(G(S^m;A,B,n);\kk), \quad i=0, 1, \dots , n;\]
the Poincar\'e polynomial of $G(S^m;A,B,n)$ equals $1+t^{m-1}+t^{2(m-1)} +\dots +t^{n(m-1)}$.
For $m\ge 3$ odd, the multiplication is given by 
\begin{eqnarray}
\sigma_i\sigma_j \, =\, \left \{
\begin{array}{l}
\displaystyle{\frac{(i+j)!}{i!\cdot j!}}\cdot \sigma_{i+j},\quad\mbox{if}\quad i+j\le n,\\ \\
0, \quad \mbox{if}\quad  i+j>n. 
\end{array}
\right .
\end{eqnarray}
For $m\ge 2$ even, the multiplication is given by
\begin{eqnarray}
\quad \quad \sigma_i\sigma_j \, =\, \left \{
\begin{array}{l}
\displaystyle{\frac{[(i+j)/2]!}{[i/2]!\cdot [j/2]!}}\cdot \sigma_{i+j},\quad\mbox{if $i+j\le n$ and $i$ or $j$ is even,}\\ \\
0, \quad \mbox{if either $i+j>n$, or both $i$ and $j$ are odd.}
\end{array}
\right.
\end{eqnarray}
\end{theorem}

Recall that $[x]$ denotes the integer part of $x$.

The proof of Theorems \ref{openstring4} is given in the following sections.

\begin{remark}
Choosing an arbitrary point $C\in S^m$, where $C\ne A,$ and $C\ne B$ we obtain
an inclusion $\phi: G(S^m-C;A,B,n)\to G(S^m;A,B,n)$; here we may identify $S^m-C$ with $\R^m$.
From the proof of Theorem \ref{openstring4} it will be clear that the induced homomorphism
$$\phi^\ast : H^\ast(G(S^m;A,B,n);\kk) \to H^\ast(G(S^m-C;A,B,n);\kk)= H^\ast(G(\R^m;A,B,n);\kk),$$
is injective and its image may easily be described.
For example, for $m$ odd, $\phi^\ast$ maps each generator $\sigma_r$ to the degree $r$ symmetric function of classes
$s_i$:
\begin{eqnarray}
\phi^\ast(\sigma_r) = \sum_{0\le i_1<i_2<\dots i_r\le n} s_{i_1}s_{i_2}\dots s_{i_r},\quad r=1, \dots, n,
\end{eqnarray}
where $s_0, \dots, s_n\in H^{m-1}(G(S^m-C;A,B,n);\kk) = H^{m-1}(G(\R^m;A,B,n);\kk)$
are the generators given by Proposition \ref{cor13}.

For even $m$ the classes $\phi^\ast(\sigma_r)$ may also be described. Such description may easily 
be extracted from the proof
of Theorem \ref{openstring4} given below. 
For instance, 
$$\phi^\ast(\sigma_1) = s_n -s_{n-1}+\dots +(-1)^n s_0\quad \mbox{and }\quad 
\phi^\ast(\sigma_2) =\, -\sum_{0\le i<j \le n} (-1)^{i+j} s_i s_j,$$
as follows from formulae (\ref{set}). More generally,
\begin{eqnarray}\label{general}
\qquad \phi^\ast(\sigma_{r}) =\, (-1)^{[r/2]+nr}\cdot \sum_{0\le i_1<i_2<\dots <i_{r}\le n}
(-1)^{i_1+i_2+\dots +i_{r}} s_{i_1}s_{i_2}\dots s_{i_{r}},
\end{eqnarray}
where the sum is taken over all increasing sequences $0\le i_1<i_2<\dots <i_{r}\le n$. 
This follows using 
$$\beta_j\beta_{j+1} =\beta_j s_{j+1}=\sum_{i=0}^j (-1)^{i+j}s_is_{j+1}$$ 
from 
our definition (\ref{even}).
\end{remark}

\section{\bf Proof of Theorem \ref{openstring4} for $m$ odd.}

Theorem \ref{openstring1}
gives a spectral sequence of bigraded algebras, which converges to the cohomology
algebra $H^\ast(G(S^m;A,B,n);\kk)$; 
the initial term $E_2=E_m$ has generators $u_1, \dots, u_n$, having bidegree $(m,0)$, which satisfy
$$
 u_iu_j =-u_ju_i,\quad u_i^2=0,
$$
and also generators $s_0, s_1, \dots, s_n$, having bidegree $(0,m-1)$, which satisfy
\begin{eqnarray}
\begin{array}{l}
s_is_j =s_js_i, \quad s_i^2=0,\\
s_iu_j=u_js_i,\\
u_1s_0=0,\\
(u_i-u_{i+1})s_i=0,\quad \mbox{for}\quad i=1, 2, \dots, n-1,\\
u_ns_n=0,\\
s_0s_1\dots s_n=0.
\end{array}
\end{eqnarray}
Here 
$u_i $ denotes $1\times \dots \times u\times 1\times \dots 1\in H^m((S^m)^{\times n};\kk)$, 
where $u$ is the fundamental class of the
sphere $S^m$, which appears on the place number $i$.

The differential $d=d_m\, :\,  E_m\to E_m$ acts by
$$
\begin{array}{l}
du_j=0,\\
ds_0=-u_1,\\
ds_i=u_i -u_{i+1}, \quad\mbox{for}\quad i=1, 2, \dots, n-1,\\
ds_n=u_n.
\end{array}
$$
Introduce new variables $v_0, v_1, \dots, v_n$:
$$
\begin{array}{l}
v_0= -u_1,\\
v_i=u_i-u_{i+1}, \quad \mbox{for}\quad i=1, \dots, n-1,\\
v_n = u_n.
\end{array}
$$

We have the following relations:
\begin{eqnarray}\label{stam}
\begin{array}{ll}
{\rm (i)} & \quad v_iv_j =-v_jv_i,\quad v_i^2=0,\\
{\rm (ii)} & \quad v_0+\dots +v_n =0,\\
{\rm (iii)} & \quad v_is_i =0,\quad \mbox{for}\quad i=0, 1, \dots, n,\\
{\rm (iv)} & \quad s_is_j=s_js_i,\quad s_iv_j=v_js_i,\\
{\rm (v)} & \quad s_i^2=0,\\
{\rm (vi)} & \quad s_0s_1\dots s_n=0,\\
{\rm (vii)} & \quad ds_i=v_i,\quad \mbox{for}\quad i=0, 1, \dots, n,\\
{\rm (viii)} & \quad dv_i=0.
\end{array}
\end{eqnarray}

Denote by $\sigma_k\in E_m$ the $k$-th symmetric function in
variables
$s_0$, $s_1,$ $\dots, s_n$, i.e.,
\[\sigma_0=1, \quad{\rm and}\quad
\sigma_k =\sum_{0\le i_1<\dots <i_k\le n} s_{i_1}s_{i_2}\dots
s_{i_k}\quad{\rm for}\quad
 k=1, 2, \dots, n.\]
It is clear (because of (v) in (\ref{stam})) that
\begin{eqnarray}\label{product}
\sigma_i\sigma_j \, =\, \left\{
\begin{array}{l}
\displaystyle
\left(\begin{array}{c} {i+j}\\ i\end{array}\right)
\sigma_{i+j}, \quad{\rm for}\quad i+j\le n,\\ \\
0,\quad{\rm for}\quad i+j>n.
\end{array}
\right.
\end{eqnarray}
It is clear that 
$$d\sigma_1=d(s_0+\dots +s_n) = v_0+\dots +v_n=0$$
and similarly
$$d(\sigma_i) = (v_0+v_1+\dots+v_n)\sigma_{i-1}=0$$
for any $i$. Hence we have found nonzero cycles $\sigma_0, \sigma_1, \dots, \sigma_n$, which 
(by the obvious geometric reasons)
cannot belong to the image of $d$. Our purpose is to show that these classes
additively generate the whole cohomology $H^\ast(E_m,d)$.

Denote by $(A_n,d)$ the graded differential algebra with generators $v_0, \dots, v_n$ and $s_0, \dots, s_n$
which satisfy relations (i), (iii), (iv), (v), (vi) among (\ref{stam}).  The differential $d:A_n\to A_n$ 
is given by formulae (vii) and (viii) of (\ref{stam}). We consider $A_n$ with the total grading, where each $s_i$ 
has degree
$m-1$ (even) and each $v_i$ has degree $m$ (odd). 

We claim that $H^j(A_n,d)=0$ for $j>0$. The proof will use induction on $n$. For $n=0$ the claim is obvious.
We have a natural inclusion $A_{n-1}\to A_n$, which identifies $A_{n-1}$ with the subalgebra of $A_n$ generated by
$s_0,\dots, s_{n-1}$ and $v_0, \dots, v_{n-1}$. We will show that the factor $A_n/A_{n-1}$ is acyclic which clearly 
implies step of induction.
Any element $a\in A_n/A_{n-1}$ can be uniquely represented in the form $a=s_nx+v_ny$, where $x, y\in A_{n-1}$.
If $da=0$ then
\[d(a)=v_nx+s_nd(x) -v_nd(y) = s_nd(x) +v_n[x-d(y)] =0\]
and hence $x=d(y)$ and $a=d(s_ny)$. The claim follows.

Introduce a new differential $\delta_n: A_n\to A_n$ of degree $m$:
$$\delta_n(x)=(\sum_{i=0}^n v_i)x.$$
Clearly, $\delta_n^2=0$ and $\delta_n d = - d\delta_n$; however $\delta_n$
does not obey the
Leibnitz rule. We claim that
\begin{eqnarray}
H^i(A_n,\delta_n)=\left\{
\begin{array}{l}
\kk,\quad {\rm if}\quad i=(n+1)(m-1),\\ \\
0, \quad {\rm otherwise}
\end{array}
\right.\label{acyc3}
\end{eqnarray}
and a nontrivial cohomology class is represented by the product 
$s_0s_1s_2\dots s_n$.
Indeed, each element of $A_n$ can be written as a sum of monomials in $s_i$,
$v_i$.
For $I\subset \{0, 1, 2, \dots n\}$,  denote by $s_I$ the product of all $s_i$
for $i\in I$. Similarly, we will label the monomials $v_{i_1}\dots v_{i_r}$ with $i_1<i_2<\dots <i_r$
as $v_J$, where $J$ is a subset
$J=\{i_1, \dots, i_r\}\subset \{0, 1, 2, \dots n\}$. Note that the product
$s_I v_J\in A_n$ is nontrivial if and only if $I$ and $J$ are disjoint
subsets of $\{0, 1, \dots, n\}$.
Note also that
\[\delta_n(s_Iv_J)\, =\,  \sum_{i\notin I\cup J} \epsilon_i s_Iv_{J\cup \{i\}},\]
where $\epsilon_i$ is $\pm 1$ depending on whether $J$ contains even or odd number of members less than $i$.
We see that application of $\delta_n$ does not change the multi-index $I$.
Hence, the complex $(A_n, \delta_n)$ splits into a direct sum over different
multi-indices $I$.
Fix a set $I$ and denote by $k$ the cardinality of the set $\{0, 1,\dots,n\} -I$.
Then the respective part of the complex $(A_n, \delta_n)$ is isomorphic to
the standard cochain complex of the simplex with $k$ vertices:

the differential of an  $r$-dimensional face (i.e., set $J$) is the sum of
$r+1$-dimensional faces
that contain the given one (sets $J \cup \{i\}$). Note that empty set $J$
is also allowed.
This complex has zero cohomology, unless $k=0$
(empty simplex), in which case the cohomology is ${\kk}$.
This exceptional case corresponds to $I = \{0,1,\dots,n\}$, and (\ref{acyc3}) follows.

Let $\I_n\subset A_n$ and $\K_n\subset A_n$ denote the image and the kernel of $\delta_n: A_n\to A_n$.
Note that $\I_n\subset \K_n$ and the factor $\K_n/\I_n$ is one-dimensional generated by 
the product $s_0s_1\dots s_n$. Hence we obtain that
\[H^j(\I_n,d) \simeq H^j(\K_n,d), \quad j\ne (n+1)(m-1),\]
and 
$H^{(n+1)(m-1)}(\I_n,d) =0.$

Since we know that $H^j(A_n,d)=0$ for $j>0$, the short exact sequence 
\[0\to \K_n\to A_n\stackrel{\delta_n}\longrightarrow \I_n\to 0\]
gives isomorphisms
\[H^{j+m-1}(\I_n,d) \simeq H^j(\K_n, d)\]
for all $j > 1$. This leads to periodicity
\[H^j(\I_n,d) \simeq H^{j+m-1}(\I_n,d), \quad \mbox{for all}\quad j\ne 1, \, \, j\ne (n+1)(m-1).\]
On the other hand it is obvious that for $1<j< 2m-1$ the cohomology $H^j(\I_n,d)$ vanishes
unless $j=m$ and for $j=m$ it is one dimensional (generated by the class $v_0+v_1 +\dots +v_n$).
This shows that 
$$
\dim H^j(\I_n,d) \, =\, \left\{
\begin{array}{l}
1, \quad \mbox{for}\quad j=i(m-1)+1, \quad 1\le i\le n+1,\\
0, \quad \mbox{otherwise.}
\end{array}
\right .
$$

Using  $H^j(A_n,d)=0$ we get 
$$
\dim H^j(A_n/\I_n,d) \, =\, \left\{
\begin{array}{l}
1, \quad \mbox{for}\quad j=i(m-1), \quad 0\le i\le n+1,\\
0, \quad \mbox{otherwise.}
\end{array}
\right .
$$

The term $E_m$ is obtained from $A_n/\I_n$ by factoring out with respect to the ideal generated by the product 
$s_0s_1\dots s_n$, 
generating the top-dimensional cohomology space $H^{(n+1)(m-1)}(A_n/\I_n,d)$. Hence 
$$
\dim H^j(E_m,d) \, =\, \left\{
\begin{array}{l}
1, \quad \mbox{for}\quad j=i(m-1), \quad 0\le i\le n,\\
0, \quad \mbox{otherwise.}
\end{array}
\right .
$$
This proves that the classes $\sigma_0, \sigma_1, \dots, \sigma_n\in H^\ast(E_m,d)$ (which were described 
in the beginning of the proof) span the cohomology.

\section{\bf Proof of Theorem \ref{openstring4} for $m$ even}
\subsection{}
Theorem \ref{openstring1}
gives a spectral sequence of bigraded algebras converging to $H^\ast(G(S^m;A,B,n);\kk)$,
with the initial term $E_2=E_m$ described below.

$E_m$ has generators $u_1, \dots, u_n$, having bidegree $(m,0)$, which satisfy
$$
u_i^2 =0,\quad u_iu_j =u_ju_i,
$$
and also generators $s_0, s_1, \dots, s_n$, having bidegree $(0,m-1)$, which satisfy
\begin{eqnarray}
\begin{array}{l}
s_is_j =-s_js_i,\quad s_i^2=0,\\
s_iu_j=u_js_i,\\
u_1s_0=0,\\
(u_i-u_{i+1})s_i=0,\quad \mbox{for}\quad i=1, 2, \dots, n-1,\\
u_ns_n=0,\\
s_0s_1\dots s_n=0.
\end{array}
\end{eqnarray}\label{relations}
Here as in the previous section,
$u_i $ denotes $1\times \dots \times u\times 1\times \dots 1\in H^m((S^m)^{\times n};\kk)$, 
where $u$ is the fundamental class of the
sphere $S^m$ and it appears on the place number $i$.

The differential $d=d_m\, :\,  E_m\to E_m$ is given by
$$
\begin{array}{l}
du_j=0,\\
ds_0=u_1,\\
ds_i=u_i +u_{i+1}, \quad\mbox{for}\quad i=1, 2, \dots, n-1,\\
ds_n=u_n.
\end{array}
$$

Our purpose is to compute the cohomology of $(E_m,d)$; from the answer which we will obtain,
it will be clear that
all further differentials $d_r, \, r>m$ vanish and thus $H^\ast(E_m^{\ast,\ast},d) = E_\infty^{\ast,\ast}$.

\subsection{} In this subsection we will describe nontrivial classes
\[\sigma_i\in E_m^{i(m-1),0}, \quad i=1, 2, \dots, n ,\]
which are cocycles $d\sigma_i=0$. They appear in the first column of the spectral sequence
and hence they cannot belong to the image of $d$. Later we will show that the cohomology classes of cocycles
$\sigma_i$
span the whole cohomology of $(E_m,d)$.

Let us denote 
\[\beta_i = s_i - s_{i-1} + \dots +(-1)^i s_0\in E^{m-1, 0}_m, \quad \mbox{for}\quad i=0, 1, \dots, n.\]
We may express $s_i$ as $\beta_i +\beta_{i-1}$ for $i\ge 1$ and $s_0=\beta_0$.
We have 
$$
\begin{array}{l}
\beta_i\beta_j =-\beta_j\beta_i,\quad \beta_i^2=0,\\
d\beta_i =u_{i+1}, \quad\mbox{for}\quad i=0, 1, \dots, n-1,\\
d\beta_n=0
\end{array}
$$
Relations (\ref{relations}) give 
\begin{eqnarray}\label{rel2}
\begin{array}{l}
\beta_0u_1 =0,\\
\beta_{i-1}u_i - \beta_iu_{i+1}+d(\beta_{i-1}\beta_i) =0,\quad \mbox{for}\quad i=1, \dots, n-1,\\
\beta_{n-1}u_n+d(\beta_{n-1}\beta_n) =0.
\end{array}
\end{eqnarray}

Now we will set 
\begin{eqnarray}\label{set}
\sigma_1 = \beta_n\in E^{m-1,0}_m, \qquad \sigma_2 =\sum_{i=0}^{n-1}\beta_i\beta_{i+1}\in E^{2(m-1),0}_m.
\end{eqnarray}
Then  $d\sigma_1 =0$ and (using (\ref{rel2})) we obtain $d\sigma_2=0$. 

For any $k \le n/2$ we will define 
\begin{eqnarray}\label{even}
\sigma_{2k} =\sum \beta_{i_1}\beta_{i_1+1}\beta_{i_2}\beta_{i_2+1}\dots \beta_{i_k}\beta_{i_k+1}\, \in \, 
E^{2k(m-1),0}_m,
\end{eqnarray}
where 
$$i_r +1 < i_{r+1},\quad 0\le i_r<n \quad \mbox{for}\quad r=1, \dots, k.$$
For  $2k+1 \le n$ we will define
$$\sigma_{2k+1} = \sigma_1 \cdot\sigma_{2k} \, \in \, E^{(2k+1)(m-1),0}_m.$$
It is clear that the classes $\sigma_{1}, \sigma_2, \dots, \sigma_n$ are nonzero and  
\begin{eqnarray}
(\sigma_2)^k = k! \cdot \sigma_{2k}\quad\mbox{and}\quad \sigma_1^2=0.
\end{eqnarray}\label{rel3}
Hence classes $\sigma_i$ satisfy the following multiplication law:
\begin{eqnarray}
\qquad \sigma_i\sigma_j \, =\, \left \{
\begin{array}{l}
0, \quad \mbox{if either $i+j>n$, or both $i$ and $j$ are odd,}\\ \\
\displaystyle{\frac{[(i+j)/2]!}{[i/2]!\cdot [j/2]!}}\cdot \sigma_{i+j},\quad\mbox{if \quad $i+j\le n$ and $i$ or $j$ is even.}
\end{array}
\right .
\end{eqnarray}

We have 
\begin{eqnarray*}
\begin{array}{ll}
d(\sigma_{2k})  &= \sum \beta_{i_1}\beta_{i_1+1}\dots \beta_{i_{k-1}}\beta_{i_{k-1}+1}d(\beta_j\beta_{j+1})=\\ \\
&= \sigma_{2(k-1)}\cdot \sum_{j=0}^{n-1}d(\beta_j\beta_{j+1}) =\sigma_{2(k-1)}\cdot d(\sigma_2)= 0.
\end{array}
\end{eqnarray*}
In the first sum $j$ runs over the set $\{0, \dots, n-1\}$ and indices $i_1, \dots, i_{k-1}$ satisfy
$0\le i_r<n$ and $i_r +1<i_{r+1}$.  
Thus we have
$$
d(\sigma_i)=0, \quad \mbox{for all}\quad i=1, \dots, n.
$$

\subsection{} Next we will show that $H^\ast (E_m,d)$ contains no nontrivial cohomology classes
except linear combinations of $\sigma_1, \dots, \sigma_{n}$. More precisely, we will show that the cohomology
$H^j(E_m,d)$ (considered with respect to the total grading)
vanishes if $j>n(m-1)$ or if $j$ is not divisible by $m-1$ and it is one-dimensional
otherwise. 

We will introduce new variables $v_j$, where $j=0, 1, \dots, n$, given by
$$
\begin{array}{l}
v_0 = u_1,\\
v_i=u_{i+1}-u_i,\quad \mbox{for}\quad i=1, 2, \dots, n-1,\\
v_n = - u_n.
\end{array}
$$
The new variables commute $v_iv_j=v_jv_i$ and satisfy:
\begin{eqnarray}\label{newrel}
\begin{array}{lll}
&{\rm (i)} & v_0^2 =0,\\
&{\rm (ii)} &v_i^2 +2v_i(v_0+v_1+\dots + v_{i-1})=0,\quad\mbox{for}\quad i=1, 2, \dots, n,\\
&{\rm (iii)} & v_0+v_1+\dots +v_n =0,\\
&{\rm (iv)} & v_is_i=0\quad\mbox{for}\quad i=0, 1, 2, \dots, n,\\
&{\rm (v)} & s_is_j=-s_js_i,\\
&{\rm (vi)} & s_0s_1\dots s_n=0,\\
&{\rm (vii)} & ds_i = 2v_0+2v_1+\dots + 2v_{i-1} +v_i, \quad\mbox{for}\quad i=0, 1, \dots, n,\\
&{\rm (viii)} & dv_i=0.
\end{array}
\end{eqnarray}

Let us denote by $(A_n,d)$ the graded differential algebra with generators $v_0, \dots, v_n$ and $s_0, \dots, s_n$
satisfying relations (i), (ii), (iv), (v). Thus, we simply ignore relations (iii) and (vi). 

The differential $d: A_n\to A_n$ is given by 
formulae (vii) and (viii). Note that the ideal generated by the relations (i), (ii), (iv), (v) is invariant under 
the differential $d$; for example, $d(v_is_i) = v_i(2v_0+\dots +2v_{i-1}+v_i)$ belongs to the ideal because of 
relation
(ii). Thus $d: A_n\to A_n$ is well defined.

\begin{lemma} $H^j(A_n,d)=0$ for all $j>0$.\label{lemma17}
\end{lemma}
\begin{proof} Using relations (i), (ii), (iv), (v) we see that
the additive basis of $A_n$ is given by monomials of the form $v_Is_J$, where 
\[I, \, J\subset \{0, 1, \dots, n\}, \quad I\cap J =\emptyset,\]
are disjoint multi-indices. Hence it is clear that for $j<n$ the differential algebra $A_j$ can be 
embedded into $A_n$; in fact $A_j$ may be identified with the subalgebra generated by
$s_0, \dots, s_j$ and $v_0, \dots, v_j$.

The factor $A_j/A_{j-1}$ has a very simple structure. Each element $a\in A_j/A_{j-1}$ has a unique
representation in the form $a = s_jx +v_jy$, where $x, y\in A_{j-1}$. From formula (v) we obtain that the
differential of $A_j/A_{j-1}$ acts as follows $da=v_jx - s_jd(x) +v_jd(y)$. Hence $da=0$ 
is equivalent to $x+dy=0$, which implies that  $a=d(s_jy)$. Thus we obtain that 
each factor $A_j/A_{j-1}$
is acyclic.

The statement of the Lemma now follows by induction.
\end{proof}

Consider now the homomorphism $\delta_n: A_n \to A_n$ given by multiplication by $v_0+v_1+\dots +v_n$, i.e.
\[\delta_n(x) = (v_0+v_1+\dots +v_n)x, \quad x\in A_n.\]
Using relations (i) and (ii) one obtains
$\delta_n^2 =0,$
i.e. $\delta_n$ may be viewed as a new differential on $A_n$. 
Note that $\delta_n$ increases the total grading by $m$.

\begin{lemma}\label{lemma18}
$$
H^j(A_n,\delta_n) \, =\, \left\{
\begin{array}{l}
0, \quad\mbox{for}\quad j\ne (n+1)(m-1),\\
\kk, \quad\mbox{for}\quad j= (n+1)(m-1)
\end{array}
\right .
$$
and the product $s_0s_1\dots s_n\in A_n$ is a cocycle (with respect to $\delta_n$), representing a nontrivial
cohomology class.
\end{lemma}
\begin{proof} We will use induction on $n$. The statement is trivial when $n=0$. 
Let's assume that it is true for $n-1$. Consider the homomorphism
\[\phi: A_{n-1} \to A_n,\quad \phi(x) = v_n\cdot x, \quad x\in A_{n-1}.\]
It is clear that $\phi$ is injective and increases the total degree by $m$. 
Using relation (ii) one finds
$$
\begin{array}{l}
\delta_n(\phi(x)) = (v_0+\dots + v_n)v_n x = \\
= -v_n(v_0+\dots +v_{n-1})\cdot x = \\
= -\phi(\delta_{n-1}(x)).
\end{array}
$$
Hence we obtain a short exact sequence
\[ 0 \to A_{n-1} \stackrel {\phi}\longrightarrow A_n \to A_n/\phi(A_{n-1}) \to 0\]
and a long homological sequence
$$
\begin{array}{l}
\stackrel {\phi}\longrightarrow H^{j}(A_n,\delta_n) \to 
H^{j}(A_n/\phi(A_{n-1}),\delta_n) 
 \stackrel {\kappa}\longrightarrow H^{j}(A_{n-1}, \delta_{n-1})\to \dots
\end{array}
$$
We will show that the connecting homomorphism
\[\kappa: H^{j}(A_n/\phi(A_{n-1}),\delta_n) \to H^{j}(A_{n-1}, \delta_{n-1})\]
is an isomorphism for all $j\ne (n+1)(m-1)$ and its is an epimorphism with 
one-dimensional kernel for $j = (n+1)(m-1)$. This clearly implies the statement of the Lemma.

Any element $a\in A_n/\phi(A_{n-1})$ has a unique representation in the form
\[a = x + s_ny, \quad x, y\in A_{n-1}.\]
Then $\delta_{n-1}(a)\in A_n/\phi(A_{n-1})$ equals $\delta_{n}(x) -s_n\delta_{n-1}(y)$ and hence
we obtain
\[H^j(A_n/\phi(A_{n-1}),\delta_n) \simeq H^j(A_{n-1},\delta_{n-1}) \oplus H^{j-m+1}(A_{n-1},\delta_{n-1}),\]
where the first summand corresponds to the class of $x$ and the second summand corresponds to the class of $y$.

Suppose that $a$ is a cycle of the relative complex 
$A_{n}/\phi(A_{n-1})$. In order to calculate $\kappa(a)$, the image under the connecting homomorphism,
 we have to view $a=x+s_ny$ as a chain in $A_n$ and compute $\delta_n(a)\in A_n$. We obtain $\delta_n(a) =\phi(a)$,
which shows that $\kappa$ is always an epimorphism and it is an isomorphism if and only if $ H^{j-m+1}(A_{n-1},\delta_{n-1})=0$;
by our induction hypothesis it holds if $j-m+1\ne n(m-1)$.

This completes the proof.
\end{proof}

Let $\I_n\subset A_n$ and $\K_n\subset A_n$ denote the image and the kernel of $\delta_n: A_n\to A_n$.
Note that $\I_n\subset \K_n$ and by Lemma \ref{lemma18} the factor $\K_n/\I_n$ is one-dimensional generated by 
the product $s_0s_1\dots s_n$. Hence we obtain that
\[H^j(\I_n,d) \simeq H^j(\K_n,d), \quad j\ne (n+1)(m-1),\]
and 
$H^{(n+1)(m-1)}(\I_n,d) =0.$

From Lemma \ref{lemma17} and short exact sequence 
\[0\to \K_n\to A_n\stackrel{\delta_n}\longrightarrow \I_n\to 0\]
we obtain isomorphisms
\[H^{j+m-1}(\I_n,d) \simeq H^j(\K_n, d)\]
for all $j > 1$. This gives periodicity
\[H^j(\I_n,d) \simeq H^{j+m-1}(\I_n,d), \quad \mbox{for all}\quad j\ne 1, \, \, j\ne (n+1)(m-1).\]
On the other hand it is obvious that for $1<j< 2m-1$ the cohomology $H^j(\I_n,d)$ vanishes
unless $j=m$ and for $j=m$ it is one dimensional (generated by the class $v_0+v_1 +\dots +v_n$).
This shows that 
$$
\dim H^j(\I_n,d) \, =\, \left\{
\begin{array}{l}
1, \quad \mbox{for}\quad j=i(m-1)+1, \quad 1\le i\le n+1,\\
0, \quad \mbox{otherwise.}
\end{array}
\right .
$$
Using Lemma \ref{lemma17} we get 
$$
\dim H^j(A_n/\I_n,d) \, =\, \left\{
\begin{array}{l}
1, \quad \mbox{for}\quad j=i(m-1), \quad 0\le i\le n+1,\\
0, \quad \mbox{otherwise.}
\end{array}
\right .
$$

\subsection{End of the proof of Theorem \ref{openstring4} for $m$ even}
The differential algebra $(E_m,d)$ is obtained from $(A_n,d)$ by adding relations (iii) and (vi) of (\ref{newrel});
therefore, $(E_m,d)$ is obtained from $(A_n/\I_n, d)$ by adding relation (vi) of (\ref{newrel}). We know that algebra
$H^\ast(A_n/\I_n,d)$ is generated by $\sigma_1, \dots, \sigma_n$, where $\deg(\sigma_i) =i(m-1)$.
It is clear that the product $s_0s_1\dots s_n$ is a nontrivial cycle of $A_n/\I_n$ having degree $(n+1)(m-1)$.
Comparing all this information we conclude that the classes $\sigma_1, \dots, \sigma_n$ form an additive
basis of $H^\ast(E_m,d)=E_{m+1}$. All further differentials $d_r$ with $r>m$ vanish. 

This clearly concludes the proof of Theorem
\ref{openstring4} for $m$ even.

\section{\bf Proof of Theorem \ref{thm1}}\label{proof}

For $m>1$
Theorem \ref{thm1} follows from Corollary \ref{cor1} and Theorem \ref{openstring4}.
If $m>1$ is odd we obtain a nonzero power
\[\sigma_1^n \, =\, n!\cdot \sigma_n \ne 0\in H^{n(m-1)}(G(S^m;A,B,n);\kk),\]
where $\kk$ is a field of characteristic zero. Hence the cohomological cup-length of $G(S^m,n)$ with $\kk$
coefficients is at least
$n$ and hence the Lusternik-Schnirelman category of $G(S^m,n)$ is at least $n+1$. 

For $m$ even we use Theorem \ref{openstring4}. It gives a nontrivial cup-product 
\[\left\{
\begin{array}{l}
\sigma_2^{n/2} =\, (n/2)!\cdot \sigma_n, \quad \mbox{if $n$ is even},\\
\sigma_1\sigma_2^{(n-1)/2} = [n/2]! \cdot \sigma_n, \quad \mbox{if $n$ is odd}.
\end{array}
\right.
\]
Hence we obtain that the Lusternik - Schnirelman category of $G(S^m;A,B,n)$ is at least
$[(n+1)/2]+1$.

In case $m=1$ we may use a direct argument. We may identify $S^1$
with the unit circle on the complex plane $\C$. Then a configuration $(x_1, \dots, x_n)\in G(S^1;A,B,n)$
(where we assume that $A\ne B$)
can be described by a point of the open $n$-dimensional unit cube $(\phi_1, \dots,  \phi_n)\in (0,1)^n$, such that
\[x_1=A\exp(2\pi i\phi_1), \quad\mbox{and}\quad x_j=x_{j-1}\exp(2\pi i \phi_j)\quad \mbox{for}\quad j=2, \dots, n.\]
If $\psi\in (0,1)$ is such that $B=A\exp(2\pi i\psi)$ then a point $(\phi_1, \dots, \phi_n)\in (0,1)^n$
corresponds to a configuration of the open string configuration space $G(S^1;A,B,n)$ if and only if
$\sum_{j=1}^n \phi_j -\psi$
is not an integer. The hyperplanes 
\[\sum_{j=1}^n \phi_j =\psi +k, \quad \mbox{where}\quad k=0, 1, \dots, n-1,\]
divide the cube $(0,1)^n$ into $n+1$ connected components, each is convex and hence contractible.
We obtain that the configuration space
$G(S^1;A,B,n)$ has $n+1$ path connected components and each is contractible. 
This gives 
$$\cat(G(S^1;A,B,n)) = n+1$$ 
and our statement follows from Corollary \ref{cor1}. \qed

\bibliographystyle{amsalpha}

\end{document}